\documentclass[11pt]{article}
\usepackage{amsmath,amsthm,amssymb}

\newcommand{\Supp}{\mathrm{Supp}\kern1.2pt}

\newcommand{\ch}{\mathrm{ch}\kern.8pt}
\newcommand{\td}{\mathrm{td}\kern.8pt}
\newcommand{\codim}{\mathrm{codim}\kern1.2pt}
\newcommand{\depth}{\mathrm{depth}\kern.8pt}

\theoremstyle{definition}

\newcounter{cont}

\begin{document}

\title{On the existence of solutions for the Maxwell equations}
\author{Luigi Corgnier}
\date{22 april 2009}
\maketitle

\begin{abstract}
Mathematical proofs are presented concerning the existence of solutions of
the Maxwell equations with suitable boundary conditions. In particular it is
stated that the well known \lq\lq delayed potentials\rq\rq provide effective
solutions of the equations, under reasonable conditions on the sources of
the fields.
\end{abstract}


\section{Introduction}

Any advanced test on theoretical classical electromagnetism (see for example
[1], [2], [3], [4]) states that all the physical properties of the
electromagnetic fields can be mathematically deduced from the Maxwell
Equations, and performs such deduction obtaining formulae (the \lq\lq
delayed potentials\rq\rq) that allow, in principle, to evaluate the fields,
starting from the knowledge of the sources, i.e. charges and currents.

The problem is that the delayed potentials are obtained, as pointed out in
the following section, with a method which ensures that, if solutions exist
for the Maxwell Equations with suitable conditions, they must be of the
provided form. Therefore the usual deduction is able to provide uniqueness
theorems, but not existence theorems for the solutions of the Maxwell
Equations.

This is not so surprising: also in similar areas, (e. g. the Laplace
Equation for harmonic functions) the deduction of uniqueness theorems is
well easier than the deduction of existence theorems.

For the above reasons, the traditional deduction of the electromagnetism
leaves open the problem that solutions of the base equations could not
exist. The usual way to overcome this situation is to state that \lq\lq due
to the physical nature of the problem, solutions must exist\rq\rq, but
clearly this cannot be accepted from a mathematical point of view.

The present paper presents existence theorems, i.e. states that, under
suitable conditions on the sources, the Maxwell Equations admit solutions,
and moreover that the delayed potentials provide an effective solution.
Similar results could be obtained using the general theory of partial
differential equations or distribution theory, but the presented deduction
is performed without any use of such general theories, and requires no
mathematical background, besides the one required for reading any
theoretical test on electromagnetism.

\noindent The first section summarises the standard deduction of the delayed
potentials, and points outs its defects.

\noindent The second section reports the obtained results, proving that the
delayed potentials provide an effective solution.

\noindent The third section summarises the results and discusses some
extensions.

\section{Standard deduction of electromagnetism}

The Maxwell Equations for the wide space have the form:

$\mathrm{rot}\mathbf{E}+\frac{1}{c}\frac{{\partial \mathbf{H}}}{{\partial t}}%
=0\qquad \qquad \qquad \qquad \qquad \qquad \qquad \qquad $(1)

$\mathrm{rot}\mathbf{H}-\frac{1}{c}\frac{{\partial \mathbf{E}}}{{\partial t}}%
=\frac{{4\pi }}{c}\mathbf{j}$ \ \ \ \ \ \ \ \ \ \ \ \ \ \ \ \ \ \ \ \ \ \ \
\ \ \ \ \ \ \ \ \ \ \ \ \ \ \ \ \ \ \ \ \ \ \ \ \ \ \ \ \ (2)

$\mathrm{div}\mathbf{H}=0\qquad \qquad \qquad \qquad \qquad \qquad \qquad
\qquad \qquad \ \ \ \ \ $(3)

$\mathrm{div}\mathbf{E}=4\pi \rho $ \ \ \ \ \ \ \ \ \ \ \ \ \ \ \ \ \ \ \ \
\ \ \ \ \ \ \ \ \ \ \ \ \ \ \ \ \ \ \ \ \ \ \ \ \ \ \ \ \ \ \ \ \ \ \ \ \ \
\ \ \ \ \ (4)

\noindent where:

\begin{itemize}
\item $\mathbf{E}$ is the electric field

\item $\mathbf{H}$ is the magnetic field

\item $c$ is the light speed in vacuum

\item $\mathbf{j}$ is the current density

\item $\rho$ is the charge density
\end{itemize}

The last two quantities, named ``sources of the fields''), are not
independent, since they must satisfy the continuity equation

$\mathrm{div}\mathbf{j}+\frac{{\partial \rho }}{{\partial t}}=0$ \ \ \ \ \ \
\ \ \ \ \ \ \ \ \ \ \ \ \ \ \ \ \ \ \ \ \ \ \ \ \ \ \ \ \ \ \ \ \ \ \ \ \ \
\ \ \ \ \ \ \ \ \ \ \ \ \ \ \ \ \ \ (5)

\noindent whose physical meaning is the conservation of the total charge.
The main problem of the electrodynamics is the following: given $\mathbf{j}$
and $\rho $, as functions of space and time satisfying (5), find the
electric and magnetic fields. The standard approach to this problem is
summarised in the following subsections.

\subsection{Definition of the potentials}

It is trivial to prove that, given any vector function $\mathbf{A}$, the
vector function $\mathbf{H}=\mathrm{rot}\mathbf{A}$ satisfies (3). Starting
from this point, the standard developments of electrodynamics perform the
claim that, in order to satisfy (3), the magnetic field $\mathbf{H}$ must be
given by

$\mathbf{H}=\mathrm{rot}\mathbf{A}$ \ \ \ \ \ \ \ \ \ \ \ \ \ \ \ \ \ \ \ \
\ \ \ \ \ \ \ \ \ \ \ \ \ \ \ \ \ \ \ \ \ \ \ \ \ \ \ \ \ \ \ \ \ \ \ \ \ \
\ \ \ \ \ \ \ \ \ \ \ (6)

\bigskip

\noindent where $\mathbf{A}$ is a suitable vector function. A proof of this
fact is not difficult, but it is not reported, since this result is not
required in the following. Assuming that $\mathbf{H}$ is given by (6) and
substituting into (1), the following is obtained:

$\mathrm{rot}(\mathbf{E}+\frac{1}{c}\frac{{\partial \mathbf{A}}}{{\partial t}%
})=0$ \ \ \ \ \ \ \ \ \ \ \ \ \ \ \ \ \ \ \ \ \ \ \ \ \ \ \ \ \ \ \ \ \ \ \
\ \ \ \ \ \ \ \ \ \ \ \ \ \ \ \ \ \ \ \ \ (7)

\bigskip

\noindent whose consequence is that the quantity in parenthesis is the
gradient of some scalar function $-\varphi $. Therefore the following
formula is obtained:

$\mathbf{E}=-\mathrm{grad}(\varphi )-\frac{1}{c}\frac{{\partial \mathbf{A}}}{%
{\partial t}}$

Having introduced the scalar potential $\varphi $ and the vector potential $%
\mathbf{A}$, the research of the fields is restricted to those given by (6)
and (7). Consequently (1) and (3) are automatically satisfied, while by
substituting into (2) and (4) and performing some calculations the following
two equations are obtained:

\bigskip

$\Delta A-\frac{1}{{c^{2}}}\frac{{\partial ^{2}\mathbf{A}}}{{\partial t^{2}}}%
-\mathrm{grad}(\mathrm{div}\mathbf{A}+\frac{1}{c}\frac{{\partial \mathbf{A}}%
}{{\partial t}})=-\frac{{4\pi }}{c}\mathbf{j}$

$\Delta \phi -\frac{1}{{c^{2}}}\frac{{\partial ^{2}\phi }}{{\partial t^{2}}}+%
\frac{1}{c}\frac{\partial }{{\partial t}}(\mathrm{div}A+\frac{1}{c}\frac{{%
\partial \phi }}{{\partial t}})=-4\pi \rho $

Normally the last equations are simplified by assuming the Gauge condition:

$\mathrm{div}\mathbf{A}+\frac{1}{c}\frac{{\partial \varphi }}{{\partial t}}%
=0 $ \ \ \ \ \ \ \ \ \ \ \ \ \ \ \ \ \ \ \ \ \ \ \ \ \ \ \ \ \ \ \ \ \ \ \ \
\ \ \ \ \ \ \ \ \ \ \ \ \ \ \ \ \ \ \ \ \ \ \ \ (8)

\bigskip

and this assumption is heuristically justified by observing that $\varphi $
and $\mathbf{A}$ are not uniquely defined by (6) and (7). Summarising, the
restriction is done to search fields given by (6) and (7), where $\mathbf{A}$
and $\varphi $ satisfy (10). With the above assumptions, the equations
satisfied by the potentials reduce to the form

$\Delta \mathbf{A}-\frac{1}{{c^{2}}}\frac{{\partial ^{2}\mathbf{A}}}{{%
\partial t^{2}}}=-\frac{{4\pi }}{c}\mathbf{j}$ \ \ \ \ \ \ \ \ \ \ \ \ \ \ \
\ \ \ \ \ \ \ \ \ \ \ \ \ \ \ \ \ \ \ \ \ \ \ \ \ \ \ \ \ \ \ \ \ \ \ \ \ \
\ (9)

$\Delta \varphi -\frac{1}{{c^{2}}}\frac{{\partial ^{2}\varphi }}{{\partial
t^{2}}}=-4\pi \rho $ \ \ \ \ \ \ \ \ \ \ \ \ \ \ \ \ \ \ \ \ \ \ \ \ \ \ \ \
\ \ \ \ \ \ \ \ \ \ \ \ \ \ \ \ \ \ \ \ \ \ \ \ \ \ (10)

If solutions of the above equations, subjected to (8), can be found, the
conclusion can be stated that the Maxwell equations have at least one
solution, given by (6) and (7).

It is readily seen that (9) and (10) are a system of 4 independent scalar
equations; moreover they are of the similar form

$\Delta F-\frac{1}{{c^{2}}}\frac{{\partial ^{2}F}}{{\partial t^{2}}}=-G$ \ \
\ \ \ \ \ \ \ \ \ \ \ \ \ \ \ \ \ \ \ \ \ \ \ \ \ \ \ \ \ \ \ \ \ \ \ \ \ \
\ \ \ \ \ \ \ \ \ \ \ \ \ \ \ \ \ \ (11)

\subsection{Elimination of the time dependence}

Equation (11) is simplified by taking its temporal Fourier transform
(accepting that the involved functions can be transformed), and by assuming
the possibility to exchange integrations with differentiations. By standard
calculations the following equation is obtained:

$\Delta f+k^{2}f=-g$\noindent\ \ \ \ \ \ \ \ \ \ \ \ \ \ \ \ \ \ \ \ \ \ \ \
\ \ \ \ \ \ \ \ \ \ \ \ \ \ \ \ \ \ \ \ \ \ \ \ \ \ \ \ \ \ \ \ \ \ \ \ \ \
\ \ \ \ \ (12)

\bigskip

where

\begin{itemize}
\item the variable replacing $t$ in the transformed domain is named $\omega$
(it has a frequency meaning)

\item $k={\omega /c}$ , where $c$ is the light speed in vacuum

\item $f(\omega)$ and $g(\omega)$ are the Fourier transforms of $F$ and $G$,
respectively.
\end{itemize}

In (12) $k$ has simply the function of a constant parameter.

\subsection{Green Theorem}

Equation (12) is treated with the well known Green method, as is done for
the Laplace Equation for the harmonic functions. As in this case, it is
clear that the Green method provides an uniqueness theorem for the solution
of (14) with assigned values on a boundary, but does not provide any
existence theorem for such solution. The Green functions is chosen as $\frac{%
{e^{-ikd}}}{d}$, where $d=\sqrt{(x-x_{0})^{2}+(y-y_{0})^{2}+(z-z_{0})^{2}}$,
and $x_{0}$, $y_{0}$, $z_{0}$ are the coordinates of an arbitrary point $%
P_{0}$ (it is easily verified that this function is everywhere regular and
satisfies (12) except at $P_{0}$, therefore it is an appropriate Green
function).

In such a way the Green Theorem is obtained in the form

$f(P_{0})=-\frac{1}{{4\pi }}\int\limits_{D}{\frac{{ge^{-ikd}}}{d}}dV+\frac{1%
}{{4\pi }}\int\limits_{S}{(f\frac{{\partial \frac{{e^{-ikd}}}{d}}}{{\partial
n}}}-\frac{{e^{-ikd}}}{d}\frac{{\partial f}}{{\partial n}})dS$

\bigskip

\noindent and provides the solution of (12) at an arbitrary point interior
to a spatial domain $D$ with boundary $S$, using the known term $g$ and the
value of the solution and of its normal derivative taken on the boundary, in
total analogy with the Green formula for the Laplace equation.

A similar formula for the solution of (11) is obtained by applying the
inverse Fourier transform, and assuming interchangeability between
differentiation and integration:

$F(P_{0})=-\frac{1}{{4\pi }}\int\limits_{D}{\frac{{G(t-\frac{d}{c})}}{d}}dV+%
\frac{1}{{4\pi }}\int\limits_{S}{(\frac{1}{d}\frac{{\partial F}}{{\partial n}%
}}-\frac{1}{d}(\frac{1}{c}\frac{{\partial F}}{{\partial t}}+\frac{F}{d}\frac{%
{\partial d}}{{\partial n}}))|_{t-\frac{d}{c}}dS$

Finally, a behavior at infinity for $F$ is assumed that ensures the
vanishing of the second integral when the domain $G$ tends to infinity, and
the following is obtained:

$F(P_{0})=-\frac{1}{{4\pi }}\int\limits_{{}}{\frac{{G(t-\frac{d}{c})}}{d}}dV$
\ \ \ \ \ \ \ \ \ \ \ \ \ \ \ \ \ \ \ \ \ \ \ \ \ \ \ \ \ \ \ \ \ \ \ \ \ \
\ \ \ \ \ \ \ \ \ \ \ \ \ \ \ (13)

\bigskip

\noindent where presently the integral is extended to the whole space. (13)
provides an explicit solution of (11), if the assumption is made that a
solution exists satisfying all the outlined conditions. Coming back to the
equations for the potentials, (13) provides:

$\mathrm{A}(P_{0})=-\frac{1}{c}\int\limits_{{}}{\frac{{\mathbf{j}(x,y,z,t-%
\frac{d}{c})}}{d}}dV$ \ \ \ \ \ \ \ \ \ \ \ \ \ \ \ \ \ \ \ \ \ \ \ \ \ \ \
\ \ \ \ \ \ \ \ \ \ \ \ \ \ \ \ \ \ \ \ \ \ \ (14)

$\varphi (P_{0})=-\int\limits_{{}}{\frac{{\rho (x,y,z,t-\frac{d}{c})}}{d}}dV$
\ \ \ \ \ \ \ \ \ \ \ \ \ \ \ \ \ \ \ \ \ \ \ \ \ \ \ \ \ \ \ \ \ \ \ \ \ \
\ \ \ \ \ \ \ \ \ \ \ \ \ \ \ (15)

\bigskip

The above formulae are known as ``delayed potentials'', and are believed to
provide, with (6) and (7), a general solution to the electrodynamics problem.

\subsection{Criticism}

As pointed out by the above review of the standard deduction, the fields
obtained applying (6) and (7) to (14) and(15) provide a solution to the
Maxwell equation assuming that the following conditions are verified:

\begin{itemize}
\item Existence of a solution which can be derived by the potentials
according to (6) and (7) (this would not be difficult to prove, if the
assumption is made of existence of a solution of (1),(2), (3), (4))

\item Existence of a solution with suitable conditions at infinity

\item As a minor point, validity of various hypothesis concerning the
possibility of exchanging integrals (singular and with infinite domain) and
derivatives

\item Existence of the integrals appearing in (14) and (15)

\item Validity, for the potentials given by (14) and (15), of the condition
(8), with is necessary to conclude that (9) and (10) provide solutions for
the potentials.
\end{itemize}

In conclusion, the summarised deduction provides only an heuristic feeling
to have found possible solutions of the Maxwell equations.

\section{Solubility Theorems}

In order to overcome the summarised limits, the chosen strategy does not
attempt to verify the validity of the various conditions pointed out at the
end of the preceding section, but attacks directly the problem to prove
that, under suitable conditions on the field sources $\rho $ and $\mathbf{j}$%
, (14) and (15) provide functions that satisfy (8), (9) and (10), from which
it is straightforward to prove that the fields obtained by (6) and (7)
satisfy the Maxwell equations (1), (2), (3), (4).

The assumed conditions on $\rho$ and $\mathbf{j}$ are the following:

\begin{itemize}
\item Regularity, i.e. existence and continuity up to the first derivatives

\item Concerning the space variables, vanishing outside of a finite regular
domain $D$, where \lq\lq regular \rq\rq means that its boundary admits
everywhere a tangent plane

\item Concerning the time dependence, sinusoidal behaviour, or existence of
the Fourier transform, vanishing outside of a finite interval for the
transformed variable $\omega$ (in practice, finitely extended frequency
spectrum)
\end{itemize}

These conditions, in particular the second and the third, are surely
redundant, even if acceptable from a physical point of view. The deduction
will show some possibility of extensions. The proofs are reported for the
following cases, of increasing complexity:

\begin{itemize}
\item Electrostatic case

\item Magnetostatic case

\item Monochromatic case

\item General case.
\end{itemize}

\subsection{Electrostatic case}

This is the case in which all the charges are fixed, i.e. the density $\rho $
does not depend on the time, and the current density $\mathbf{j}$ vanishes
(note that condition (5) is trivially satisfied). In this case the proposed
solutions (16) and (17) assume the form

$\mathbf{A}(P_{0})=0$ \ \ \ \ \ \ \ \ \ \ \ \ \ \ \ \ \ \ \ \ \ \ \ \ \ \ \
\ \ \ \ \ \ \ \ \ \ \ \ \ \ \ \ \ \ \ \ \ \ \ \ \ \ \ \ \ \ \ \ \ \ \ \ \ \
\ \ \ \ \ \ \ \ \ \ \ \ \ \ \ \ (14)

$\varphi (P_{0})=-\int\limits_{D}{\frac{{\rho (x,y,z)}}{d}}dV$ \ \ \ \ \ \ \
\ \ \ \ \ \ \ \ \ \ \ \ \ \ \ \ \ \ \ \ \ \ \ \ \ \ \ \ \ \ \ \ \ \ \ \ \ \
\ \ \ \ \ \ \ \ \ \ \ \ \ \ \ \ \ \ \ (15)

\bigskip

\noindent where $D$ is the finite domain in which $\rho $ does not vanish,
while the equations to be verified ((8), (9) and (10)) reduce to

$\Delta \varphi =-4\pi \rho $ \ \ \ \ \ \ \ \ \ \ \ \ \ \ \ \ \ \ \ \ \ \ \
\ \ \ \ \ \ \ \ \ \ \ \ \ \ \ \ \ \ \ \ \ \ \ \ \ \ \ \ \ \ \ \ \ \ \ \ \ \
\ \ \ \ \ \ \ \ \ \ \ \ \ \ \ \ \ \ \ \ (16)

If a point $P_{0}$ is external to $D$ , the integral appearing in (15) is
not singular and can be differentiated under the sign, providing trivially $%
\Delta \varphi =-4\pi \rho $, as required. \noindent If the point $P_{0}$ is
internal to $D$, a transformation to polar coordinates $u$, $\alpha $, $%
\beta $ centered on $P_{0}$ is performed, giving

$\varphi (P_{0})=-\int_{0}^{2\pi }{d\beta }\int_{-{\pi /2}}^{{\pi /2}}{\sin
^{2}\alpha }\int_{0}^{g(P_{0},\alpha ,\beta )}{\rho (P_{0}+\mathbf{v})udu}$
\ \ \ \ \ \ \ \ \ \ \ \ \ \ \ \ \ \ \ (17)

\bigskip

\noindent in which the vector $\mathbf{v}=(u\sin \alpha \cos \beta ,u\sin
\alpha \sin \beta ,u\cos \beta )$ has been introduced, and $g(P_{0},\alpha
,\beta )$ is the distance between $P_{0}$ and the point on the boundary of $%
D $ with polar angles $\alpha $, $\beta $ , which is a regular function.

In (17) every singularity has disappeared, and this ensures the existence
and regularity of $\varphi (P_{0})$, and the existence of the vector field
given by

\bigskip\ $\mathbf{E}=-\mathrm{grad}(\varphi )$

The problem is if $\mathbf{E}$ can be calculated by taking the derivatives
respect to $P_{0}$ under the sign of integral appearing in (19). The
response is certainly "yes" if the integrals obtained by formal derivation
under the sign are uniformly convergent. Performing the formal derivation,
one obtains

$\frac{{\partial \varphi }}{{\partial x_{0}}}=^{?}-\int\limits_{D}{\rho 
\frac{{x_{0}-x}}{{d^{3}}}}dV$

\noindent where the interrogation mark means that the result must be
justified. Performing the same transformation to polar coordinates just
used, the last integral transforms to

$\varphi (P_{0})=^{?}-\int_{0}^{2\pi }{d\beta }\int_{-{\pi /2}}^{{\pi /2}}{%
\sin ^{2}\alpha }\int_{0}^{g(P_{0},\alpha ,\beta )}{\rho (P_{0}+\mathbf{v})du%
}$

\bigskip

\noindent which is without singularities. Therefore it represents a
continuous function of $P_{0}$, and the uniform convergence is ensured,
justifying the formal derivation. It is important to note that a similar
reasoning could not be applied to the second derivatives, vanishing any
attempt to verify (20) with a direct calculation. Having justified first
order derivatives under the sign of $\varphi (P_{0})$, the following result
is obtained

$\mathbf{E}(P_{0})=\int\limits_{D}{\rho (P)\frac{{\mathbf{r}_{0}-\mathbf{r}}%
}{{|\mathbf{r}_{0}-\mathbf{r}|^{3}}}}dV$ \ \ \ \ \ \ \ \ \ \ \ \ \ \ \ \ \ \
\ \ \ \ \ \ \ \ \ \ \ \ \ \ \ \ \ \ \ \ \ \ \ \ \ \ \ \ \ \ \ \ \ \ \ \ \ \
\ \ \ \ \ \ \ \ \ \ \ \ (18)

\bigskip

\noindent with uniform convergence of the integral at its unique singular
point $P_{0}$. Now let $V$ be an arbitrary regular domain internal to $D$,
with boundary $S$. The aim is to calculate the flux of $\mathbf{E}$ out of $%
V $. With some straightforward calculations, performing an integration
exchange justified by the uniform convergence already established, one finds:

$\bigskip $

$\int\limits_{S}{\mathbf{E}.\mathbf{n}dS_{0}=\int\limits_{V}{\rho
(P)dV\int\limits_{S}{\frac{{(\mathbf{r}_{0}-\mathbf{r}).\mathbf{n}}}{{|%
\mathbf{r}_{0}-\mathbf{r}|^{3}}}}dS_{0}}}+\int\limits_{D-V}{\rho
(P)dV\int\limits_{S}{\frac{{(\mathbf{r}_{0}-\mathbf{r}).\mathbf{n}}}{{|%
\mathbf{r}_{0}-\mathbf{r}|^{3}}}}dS_{0}}$ \ \ \ \ \ \ \ \ \ \ (19)

\bigskip

\noindent where \textbf{n} is the vector of modulus 1 normal to $S$. In (19)
the second integral does not contain singularities, since $P_{0}$ is on $S$,
while $P$ is external to $S$. Since a trivial calculation gives

$\bigskip $

$\mathrm{div}(\frac{{\mathbf{r}_{0}-\mathbf{r}}}{{|\mathbf{r}_{0}-\mathbf{r}%
|^{3}}})=0$

\bigskip

\noindent when $P_{0}\neq P$, a standard application of Gauss Theorem proves
that the second integral vanishes. By the same reasoning, in the first
integral the boundary $S$ can be modified to a sphere centered on $P$, and
totally included in $D$, without modifying its value. In this conditions, a
direct calculation shows that

\bigskip

$\int\limits_{S}{\mathbf{E}.\mathbf{n}dS_{0}=4\pi \int\limits_{V}{\rho (P)dV}%
}$

\bigskip

\noindent and using again Gauss Theorem

$\bigskip $

$\int\limits_{V}{\mathrm{div}\mathbf{E}dV=4\pi \int\limits_{V}{\rho (P)dV}}$

Finally, keeping into account that the domain $V$ is arbitrary and that the
functions under integral are continuous, one obtains

\bigskip

$\mathrm{div}\mathbf{E}=4\pi \rho $

\bigskip

\noindent which, inserting the definition of $\mathbf{E}$, is equivalent to
the equation to be verified (16). The focal point in the reported prove is
that a vector function like (18) has flux zero from any closed surface not
including singularities of the integral; from this it follows that the flux
from an arbitrary surface can be calculated by modifying the surface up to a
sphere concentrated on the unique singular point, for which the calculation
is trivial. Similar methods are applied to the calculation of the residues
of an analytical function, providing Cauchy formula, whose aspect is similar
to (18). This localisation property allows also to understand how the
restriction on the finiteness of the domain $D$ could be lowered: it is
sufficient to separate $D$ into a finite part surrounding the point $P$ and
a remaining infinite part. For the finite part the reported reasoning
applies, while the remaining one gives no contribution. Obviously some
conditions must be added, which authorise the performed interchange of
integrals; a sufficient condition is the existence of the integrals
appearing in (19) as multiple integrals, by Fubini Theorem.

\subsection{Magnetostatic case}

This is the case in which both the density $\rho $ and the current density $%
\mathbf{j}$ do not depend on time. In this case it will be necessary to use
condition (5), which becomes

$\bigskip \mathrm{div}\mathbf{j}=0$ \ \ \ \ \ \ \ \ \ \ \ \ \ \ \ \ \ \ \ \
\ \ \ \ \ \ \ \ \ \ \ \ \ \ \ \ \ \ \ \ \ \ \ \ \ \ \ \ \ \ \ \ \ \ \ \ \ \
\ \ \ \ \ \ \ \ \ \ \ \ \ \ \ \ \ \ \ \ \ \ \ \ \ \ \ \ \ \ \ \ \ \ \ \ (20)

In this case the proposed solutions (14) and (15) assume the form

$\bigskip $

$\mathbf{A}(P_{0})=-\frac{1}{c}\int\limits_{D}{\frac{{\mathbf{j}(x,y,z)}}{d}}%
dV$ \ \ \ \ \ \ \ \ \ \ \ \ \ \ \ \ \ \ \ \ \ \ \ \ \ \ \ \ \ \ \ \ \ \ \ \
\ \ \ \ \ \ \ \ \ \ \ \ \ \ \ \ \ \ \ \ \ \ \ \ \ \ \ \ \ \ \ \ \ (21)

$\varphi (P_{0})=-\int\limits_{D}{\frac{{\rho (x,y,z)}}{d}}dV$ \ \ \ \ \ \ \
\ \ \ \ \ \ \ \ \ \ \ \ \ \ \ \ \ \ \ \ \ \ \ \ \ \ \ \ \ \ \ \ \ \ \ \ \ \
\ \ \ \ \ \ \ \ \ \ \ \ \ \ \ \ \ \ \ \ \ \ \ \ \ \ \ (22)

\bigskip

\noindent where $D$ is the finite domain in which $\rho $ and $\mathbf{j}$
do not vanish, while the equations to be verified ((8), (9) and (10)) reduce
to

\bigskip

$\mathrm{div}\mathbf{A}=0$ \ \ \ \ \ \ \ \ \ \ \ \ \ \ \ \ \ \ \ \ \ \ \ \ \
\ \ \ \ \ \ \ \ \ \ \ \ \ \ \ \ \ \ \ \ \ \ \ \ \ \ \ \ \ \ \ \ \ \ \ \ \ \
\ \ \ \ \ \ \ \ \ \ \ \ \ \ \ \ \ \ \ \ \ \ \ \ \ \ \ \ \ (23)

$\Delta \mathbf{A}=-\frac{{4\pi }}{c}\mathbf{j}$ \ \ \ \ \ \ \ \ \ \ \ \ \ \
\ \ \ \ \ \ \ \ \ \ \ \ \ \ \ \ \ \ \ \ \ \ \ \ \ \ \ \ \ \ \ \ \ \ \ \ \ \
\ \ \ \ \ \ \ \ \ \ \ \ \ \ \ \ \ \ \ \ \ \ \ \ \ \ \ \ \ \ \ \ \ \ \ \ \
(24)

$\Delta \varphi =-4\pi \rho $ \ \ \ \ \ \ \ \ \ \ \ \ \ \ \ \ \ \ \ \ \ \ \
\ \ \ \ \ \ \ \ \ \ \ \ \ \ \ \ \ \ \ \ \ \ \ \ \ \ \ \ \ \ \ \ \ \ \ \ \ \
\ \ \ \ \ \ \ \ \ \ \ \ \ \ \ \ \ \ \ \ \ \ \ \ \ \ \ \ (25)

(25) has been already proved in the preceding section. The same reasoning,
applied to the individual components of $\mathbf{A}$ and $\mathbf{j}$,
proves (24). Therefore all what is required is to prove (23). To this aim,
the first point is to consider two vectors $\mathbf{r}_{0}$ and $\mathbf{r}$
and to evaluate the flux of the vector $\frac{\mathbf{j}}{{|\mathbf{r}_{0}-%
\mathbf{r}|}}$ \ \noindent from the boundary $S$ of $D$. Considering firstly
the case in which $P_{0}$ is external to $D$, which avoids any singularity,
and using (20), (21), (22) and Gauss Theorem, one obtains

$\bigskip $

$\int\limits_{S}{\frac{{\mathbf{j}.\mathbf{n}}}{{|\mathbf{r}_{0}-\mathbf{r}|}%
}}dS=\int\limits_{D}{\mathrm{div}\frac{\mathbf{j}}{{|\mathbf{r}_{0}-\mathbf{r%
}|}}}dV=\int\limits_{D}{\frac{1}{{|\mathbf{r}_{0}-\mathbf{r}|}}\mathrm{div}%
\mathbf{j}}dV+\int\limits_{D}{\mathbf{j}.\mathrm{grad}\frac{1}{{|\mathbf{r}%
_{0}-\mathbf{r}|}}}dV=\int\limits_{D}{\frac{{\mathbf{j}.(\mathbf{r}_{0}-%
\mathbf{r})}}{{|\mathbf{r}_{0}-\mathbf{r}|^{3}}}}dV=\mathrm{div}_{0}\mathbf{A%
}(\mathbf{r}_{0})$

\noindent (in the last passage a derivative under the sign has been taken,
justified as in the preceding section).

Therefore in this case the proof of (23) is obtained, if the flux appearing
at the left hand side is zero. A sufficient condition for that is that at
each point of $S$ $\mathbf{j}$ is tangent to $S$, or in particular null. To
prove this, suppose that exists a point $P$ on $S$ at which $\mathbf{j}$ is
directed e.g. towards the extern. Then it should be possible to find an
neighborhood of $P$ satisfying the same condition (remember that $\mathbf{j}$
is assumed continuous), with the consequence that the flux of $\mathbf{j}$
from such neighborhood, completed to a closed surface, should be strictly
positive, in contradiction with (20).

Considering now the case in which $P_{0}$ is internal to $D$, the integral
appearing in (21) becomes singular, and the above passages are not valid.
But they can be applied excluding from the integral the contribution of a
sphere $D_{\varepsilon }$ of ray $\varepsilon $ centered on $\mathbf{r}_{0}$
and with boundary $S_{\varepsilon }$, finding

\bigskip

$0=\int\limits_{S}{\frac{{\mathbf{j}.\mathbf{n}}}{{|\mathbf{r}_{0}-\mathbf{r}%
|}}}dS=\int\limits_{S}{}+\int\limits_{S_{\varepsilon
}}{}-\int\limits_{S_{\varepsilon }}{\frac{{\mathbf{j}.\mathbf{n}}}{{|\mathbf{%
r}_{0}-\mathbf{r}|}}dS=}\int\limits_{D-D_{\varepsilon }}{\mathrm{div}}\frac{%
\mathbf{j}}{{|\mathbf{r}_{0}-\mathbf{r}|}}dV+\int\limits_{S_{\varepsilon }}{%
\frac{{\mathbf{j}.\mathbf{n}}}{{|\mathbf{r}_{0}-\mathbf{r}|}}dS}$

The first integral is not singular, and can be treated as in the first case;
for the second one the mean value theorem is applied, obtaining

\bigskip

$\int\limits_{D-D_{\varepsilon }}{}\frac{{\mathbf{j}.(\mathbf{r}_{0}-\mathbf{%
r})}}{{|\mathbf{r}_{0}-\mathbf{r}|^{3}}}dV+K^{\ast }\varepsilon =0$

\bigskip

\noindent where $K^{\ast }$ is a suitable value comprised between the
minimum and the maximum of $\mathbf{j.n}$ on $S$, and therefore bounded.
Taking the limit for $\varepsilon \rightarrow 0$, the desired result $%
\mathrm{div}_{0}\mathbf{A}(\mathbf{r}_{0})=0$ is obtained.

Finally, in the case of $P$ on the boundary of $D$, it is sufficient to use
the already obtained results and the continuity of $\mathrm{div}\mathbf{A}$,
whose prove is obtained by the method of the preceding section.

\subsection{Monocromatic case}

This is the case in which the density $\rho $ and the current density $%
\mathbf{j}$ have a sinusoidal time dependence. Using the complex exponential
notation, which simplifies some arguments, they are given by

$\bigskip $

$\rho (x,y,z,t)=\rho _{a}(x,y,z)\exp (-i\omega t)$ \ \ \ \ \ \ \ \ \ \ \ \ \
\ \ \ \ \ \ \ \ \ \ \ \ \ \ \ \ \ \ \ \ \ \ \ \ \ \ \ \ \ \ \ \ \ \ \ (26)

$\mathbf{j}(x,y,z,t)=\mathbf{j}_{a}(x,y,z)\exp (-i\omega t)$ \ \ \ \ \ \ \ \
\ \ \ \ \ \ \ \ \ \ \ \ \ \ \ \ \ \ \ \ \ \ \ \ \ \ \ \ \ \ \ \ \ \ \ \ \ \
\ \ \ \ \ (27)

\bigskip

\noindent where $\rho _{a}(x,y,z)$ and $\mathbf{j}_{a}(x,y,z)$ depend only
on the spatial variables, and $\omega $ is a fixed parameter. In this case
the proposed solutions (14) and (15) assume the form

$\bigskip $

$\mathbf{A}(P_{0})=-\frac{1}{c}\int\limits_{D}{\frac{{\mathbf{j}%
_{a}(x,y,z)\exp (i\frac{\omega }{c}d)}}{d}}dV.\exp (-i\omega t)$ \ \ \ \ \ \
\ \ \ \ \ \ \ \ \ \ \ \ \ \ \ \ \ \ \ \ \ \ \ \ \ \ \ \ \ (28)

$\varphi (P_{0})=-\int\limits_{D}{\frac{{\rho _{a}(x,y,z)\exp (i\frac{\omega 
}{c}d)}}{d}}dV.\exp (-i\omega t)$ \ \ \ \ \ \ \ \ \ \ \ \ \ \ \ \ \ \ \ \ \
\ \ \ \ \ \ \ \ \ \ \ \ \ \ \ \ \ (29)

\bigskip

\noindent where $D$ is a finite domain outside of which $\rho $ and $\mathbf{%
j}$ vanish. Posing

$\bigskip \mathbf{A}=\mathbf{A}_{a}\exp (-i\omega t)$

$\varphi =\varphi _{a}\exp (-i\omega t)$

\bigskip

equations (28) and (29) give

$\bigskip $

$\mathbf{A}_{a}(P_{0})=-\frac{1}{c}\int\limits_{D}{\frac{{\mathbf{j}%
_{a}(x,y,z)\exp (i\frac{\omega }{c}d)}}{d}}dV$ \ \ \ \ \ \ \ \ \ \ \ \ \ \ \
\ \ \ \ \ \ \ \ \ \ \ \ \ \ \ \ \ \ \ \ \ \ \ \ \ \ \ \ \ \ \ \ \ \ \ \ (30)

$\varphi _{a}(P_{0})=-\int\limits_{D}{\frac{{\rho _{a}(x,y,z)\exp (i\frac{%
\omega }{c}d)}}{d}}dV$ \ \ \ \ \ \ \ \ \ \ \ \ \ \ \ \ \ \ \ \ \ \ \ \ \ \ \
\ \ \ \ \ \ \ \ \ \ \ \ \ \ \ \ \ \ \ \ \ \ \ \ \ \ \ (31)

\bigskip

\noindent while the equations to be verified ((8), (9)and (10)) reduce to

$\bigskip $

$\mathrm{div}\mathbf{A}_{a}-ik\varphi _{a}=0$ \ \ \ \ \ \ \ \ \ \ \ \ \ \ \
\ \ \ \ \ \ \ \ \ \ \ \ \ \ \ \ \ \ \ \ \ \ \ \ \ \ \ \ \ \ \ \ \ \ \ \ \ \
\ \ \ \ \ \ \ \ \ \ \ \ \ \ \ \ \ \ \ \ \ \ \ \ (32)

$\Delta \varphi _{a}+k^{2}\varphi _{a}=-4\pi \varphi _{a}$ \ \ \ \ \ \ \ \ \
\ \ \ \ \ \ \ \ \ \ \ \ \ \ \ \ \ \ \ \ \ \ \ \ \ \ \ \ \ \ \ \ \ \ \ \ \ \
\ \ \ \ \ \ \ \ \ \ \ \ \ \ \ \ \ \ \ \ \ \ \ (33)

$\Delta \mathbf{A}_{a}+k^{2}\mathbf{A}_{a}=-\frac{{4\pi }}{c}\mathbf{j}_{a}$
\ \ \ \ \ \ \ \ \ \ \ \ \ \ \ \ \ \ \ \ \ \ \ \ \ \ \ \ \ \ \ \ \ \ \ \ \ \
\ \ \ \ \ \ \ \ \ \ \ \ \ \ \ \ \ \ \ \ \ \ \ \ \ \ \ \ \ \ \ \ \ (34)

\bigskip

\noindent where $k={\omega /c}$.Finally the continuity condition (5) becomes

$\bigskip $

$\mathrm{div}\mathbf{j}_{a}-i\omega \rho _{a}=0$ \ \ \ \ \ \ \ \ \ \ \ \ \ \
\ \ \ \ \ \ \ \ \ \ \ \ \ \ \ \ \ \ \ \ \ \ \ \ \ \ \ \ \ \ \ \ \ \ \ \ \ \
\ \ \ \ \ \ \ \ \ \ \ \ \ \ \ \ \ \ \ \ \ \ \ \ \ \ \ (35)

The first step is to verify (32). The detailed steps are not reported, since
they are identical to those used in the preceding section, with the only
difference that the starting point is the calculation of the flux of $\ \ 
\frac{{\mathbf{j}_{a}\exp (i\omega |\mathbf{r}_{0}-\mathbf{r}|)}}{{|\mathbf{r%
}_{0}-\mathbf{r}|}}$ \ \noindent which vanishes due to the tangent condition
on $\mathbf{j}_{a}$.

Coming to the proof of (33), one defines

$\bigskip $

$\mathbf{E}_{a}=-\mathrm{grad}\varphi _{a}+ik\mathbf{A}_{a}$ \ \ \ \ \ \ \ \
\ \ \ \ \ \ \ \ \ \ \ \ \ \ \ \ \ \ \ \ \ \ \ \ \ \ \ \ \ \ \ \ \ \ \ \ \ \
\ \ \ \ \ \ \ \ \ \ \ \ \ \ \ \ \ \ \ \ \ \ \ \ \ \ (36)

\bigskip

\noindent and, using (32), obtains

$\bigskip $

$\mathrm{div}_{0}\mathbf{E}_{a}=-\mathrm{div}_{0}\int\limits_{D}{\rho
_{a}(ik|\mathbf{r}_{0}-\mathbf{r}|-1)\exp (ik|\mathbf{r}_{0}-\mathbf{r}|)}%
\frac{{\mathbf{r}_{0}-\mathbf{r}}}{{|\mathbf{r}_{0}-\mathbf{r}|^{3}}}%
dV+k^{2}\varphi _{a}$

Using the fact that an integral with first order singularity in $|\mathbf{r}%
_{0}-\mathbf{r}|$ can be differentiated under the sign and developing some
calculations, the following is obtained:

$\bigskip $

$\mathrm{div}_{0}\mathbf{E}_{a}=-ik\int\limits_{D}{\rho _{a}(ik|\mathbf{r}%
_{0}-\mathbf{r}|-1)\exp (ik|\mathbf{r}_{0}-\mathbf{r}|)}\frac{1}{{|\mathbf{r}%
_{0}-\mathbf{r}|^{2}}}dV+\mathrm{div}_{0}\mathbf{C}$ \ \ \ \ \ \ \ \ (37)

\noindent having defined

$\bigskip $

$\mathbf{C}=\int\limits_{D}{\rho _{a}\exp (ik|\mathbf{r}_{0}-\mathbf{r}|)}%
\frac{{\mathbf{r}_{0}-\mathbf{r}}}{{|r_{0}-r|^{3}}}dV$

\bigskip

At this point the flux of C through a surface S contained in V is evaluated
by steps similar to those reported in the section on electrostatics. The
result is

\bigskip

$\int\limits_{V}{\mathrm{div}\mathbf{C}dV_{0}}=\int\limits_{S}{dS_{0}}%
\int\limits_{D}{\rho _{a}(P)\frac{{(\mathbf{r}_{0}-\mathbf{r}).\mathbf{n}}}{{%
|\mathbf{r}_{0}-\mathbf{r}|^{3}}}\exp (ik|\mathbf{r}_{0}-\mathbf{r}|)}%
dV=\int\limits_{V}{}+\int\limits_{D-V}{\rho _{a}(P)dV\int\limits_{S}{\frac{{(%
\mathbf{r}_{0}-\mathbf{r}).\mathbf{n}}}{{|\mathbf{r}_{0}-\mathbf{r}|^{3}}}%
\exp (ik|\mathbf{r}_{0}-\mathbf{r}|)}}dS_{0}$

\bigskip

In the above equation the second integral is regular, and its value is

\bigskip

$\int\limits_{D-V}{\rho _{a}(P)dV\int\limits_{S}{\frac{{(\mathbf{r}_{0}-%
\mathbf{r}).n}}{{|\mathbf{r}_{0}-\mathbf{r}|^{3}}}\exp (ik|\mathbf{r}_{0}-%
\mathbf{r}|)}}dS_{0}=\int\limits_{V}{dV_{0}\int\limits_{D-V}{\rho _{a}(P)%
\frac{{ik}}{{|\mathbf{r}_{0}-\mathbf{r}|^{2}}}\exp (ik|\mathbf{r}_{0}-%
\mathbf{r}|)}}dV$

\bigskip

On the contrary, in the first integral the surface is moved up to the
already used sphere of radius $\varepsilon $, and the following is obtained

\bigskip

$\int\limits_{V}{\mathrm{div}\mathbf{C}dV_{0}}=\int\limits_{V}{dV_{0}}%
\int\limits_{D-V}{\rho _{a}\frac{{ik\exp (ik|\mathbf{r}_{0}-\mathbf{r}|)}}{{|%
\mathbf{r}_{0}-\mathbf{r}|^{2}}}}dV+\int\limits_{V-V_{\varepsilon }}{dV_{0}}%
\int\limits_{V}{\rho _{a}\frac{{ik\exp (ik|\mathbf{r}_{0}-\mathbf{r}|)}}{{|%
\mathbf{r}_{0}-\mathbf{r}|^{2}}}}dV+4\pi \int\limits_{V}{\rho _{a}}dV$

\bigskip

Taking the limit $\varepsilon \rightarrow 0$ and using the fact that $V$ is
arbitrary, the final result is obtained

$\bigskip $

$\mathrm{div}\mathbf{C}=ik\int\limits_{V}{\rho _{a}\frac{{\exp (ik|\mathbf{r}%
_{0}-\mathbf{r}|)}}{{|\mathbf{r}_{0}-\mathbf{r}|^{2}}}}dV+4\pi \varphi _{a}$

\bigskip

which inserted into (37) provides the desired result (33).

The final step is the verification of (34). At a first sight, it could
appear sufficient to apply the last obtained result to the three components
of (34). However, there is a problem, due to the fact that to obtain (33)
from (31), essential use has been made of (32), which depends on the fact
that a-priori $\rho _{a}$ is not arbitrary, but tied to a certain vector $%
\mathbf{j}_{a}$ by (35). Therefore the proof is complete if the fact can be
proved that to each component $j_{a^{\prime }}$ of a vector $\mathbf{j}_{a}$
another vector $\mathbf{j}_{a}^{\ast }$ can be associated, satisfying the
equation

$\bigskip $

$\mathrm{div}\mathbf{j}_{a}^{\ast }-i\omega j_{a^{\prime }}^{{}}=0$ \ \ \ \
\ \ \ \ \ \ \ \ \ \ \ \ \ \ \ \ \ \ \ \ \ \ \ \ \ \ \ \ \ \ \ \ \ \ \ \ \ \
\ \ \ \ \ \ \ \ \ \ \ \ \ \ \ \ \ \ \ \ \ \ \ \ \ \ \ \ \ \ \ \ \ \ \ \ \ \
(38)

\bigskip

and moreover tangent, as $\mathbf{j}_{a}$, on the boundary of the domain $D$.

This is consequence of the following property:

\noindent The equation in $\mathbf{j}$

$\bigskip $

$\mathrm{div}\mathbf{j}=f(x,y,z)$ \ \ \ \ \ \ \ \ \ \ \ \ \ \ \ \ \ \ \ \ \
\ \ \ \ \ \ \ \ \ \ \ \ \ \ \ \ \ \ \ \ \ \ \ \ \ \ \ \ \ \ \ \ \ \ \ \ \ \
\ \ \ \ \ \ \ \ \ \ \ \ \ \ \ \ \ \ \ \ \ \ \ \ \ (39)

\bigskip

has solutions, even if $\mathbf{j}$ is constrained to be tangent to an
assigned surface $S$. In fact it is easy to verify that (39) (but in general
not the constraint) is satisfied by the vector

$\bigskip $

$\mathbf{j}^{\ast \ast }=(0,0,\int {f(x,y,z)dz})$

\bigskip

and therefore also by

$\bigskip $

$\mathbf{j}=\mathbf{j}^{\ast \ast }+\mathrm{rot}\mathbf{P}$

\bigskip

where $\mathbf{P}$ is an arbitrary vector. Therefore all is reduced to
choice $\mathbf{P}$ in such a way that on $S$:

$\bigskip $

$\mathrm{rot}\mathbf{P}.\mathbf{n}=-\mathbf{j}^{\ast \ast }.\mathbf{n}$

\bigskip

i.e. to prove that the equation in $P$

$\bigskip $

$\mathrm{rot}\mathbf{P}.\mathbf{n}=g(x,y,z)$ \ \ \ \ \ \ \ \ \ \ \ \ \ \ \ \
\ \ \ \ \ \ \ \ \ \ \ \ \ \ \ \ \ \ \ \ \ \ \ \ \ \ \ \ \ \ \ \ \ \ \ \ \ \
\ \ \ \ \ \ \ \ \ \ \ \ \ \ \ \ \ \ \ \ \ \ \ \ \ (40)

\bigskip

has solutions on $S$.

In fact this is possible even by imposing $P_{x}=P_{y}=0$, and reducing (40)
to the form 
\begin{equation*}
\frac{{\partial P_{z}}}{{\partial y}}n_{x}-\frac{{\partial P_{z}}}{{\partial
x}}n_{y}=g(x,y,z)
\end{equation*}%
which is a first order partial differential linear equation and can be
solved through a method which transform it into an ordinary differential
equation (see any text of Mathematical Analysis, e.g. [5]).

\subsection{General case}

This is the case in which the density $\rho $ and the current density $%
\mathbf{j}$ have arbitrary time dependence, but with the restriction to
admit a Fourier transform, producing functions vanishing outside of a finite
domain. \noindent Under these conditions, $\rho $ and $\mathbf{j}$ are of
the following form

\bigskip

$\rho (r,t)=\frac{1}{\sqrt{2\pi }}\int\limits_{B}{\rho _{a}(r,\omega )\exp
(-i\omega t)d\omega }$ \ \ \ \ \ \ \ \ \ \ \ \ \ \ \ \ \ \ \ \ \ \ \ \ \ \ \
\ \ \ \ \ \ \ \ \ \ \ \ \ \ \ \ (41)

$\mathbf{j}(r,t)=\frac{1}{\sqrt{2\pi }}\int\limits_{B}{\mathbf{j}%
_{a}(r,\omega )\exp (-i\omega t)d\omega }$ \ \ \ \ \ \ \ \ \ \ \ \ \ \ \ \ \
\ \ \ \ \ \ \ \ \ \ \ \ \ \ \ \ \ \ \ \ \ \ \ \ \ \ \ \ \ \ (42)

\bigskip

where $B$ is a finite domain. \noindent Inserting (41) and (42) into (14)
and (15) and interchanging bounded integrations, the proposed solutions take
the form

\bigskip $\mathbf{A}(P_{0},t)=-\frac{1}{{\sqrt{2\pi }c}}\int\limits_{B}{\exp
(-i\omega t)du\int\limits_{D}{\frac{{\mathbf{j}_{a}(r,\omega )\exp (i\omega }%
d/c{)}}{d}dV}}=\frac{1}{\sqrt{2\pi }}\int\limits_{B}{\mathbf{A}_{a}(r,\omega
)\exp (-i\omega t)d\omega }$ \ \ \ (43)

$\varphi (P_{0},t)=-\frac{1}{\sqrt{2\pi }}\int\limits_{B}{\exp (-i\omega
t)du\int\limits_{D}{\frac{{\rho _{a}(r,\omega )\exp (i{\omega }d/c)}}{d}dV}}=%
\frac{1}{\sqrt{2\pi }}\int\limits_{B}{\varphi _{a}(r,\omega )\exp (-i\omega
t)d\omega }$ \ \ \ \ (44){\ }

\noindent where $\mathbf{A}_{a}$ and $\varphi _{a}$ are given by (28) and
(29).

Similarly, taking derivatives under the sign of the regular range bounded
integrals (41) and (42) and using inversion properties for the Fourier
transform, the continuity equation (5) provides (32) for $\mathbf{A}_{a}$
and $\varphi _{a}$. Using the results of the preceding section, it can be
concluded that $\mathbf{A}_{a}$ and $\varphi _{a}$ satisfy (32), (33), (34).
Using this fact, one obtains, performing derivatives under bounded regular
integrals:

$\bigskip $

$\mathrm{div}\mathbf{A}+\frac{1}{c}\frac{{\partial \varphi }}{{\partial t}}=%
\frac{1}{\sqrt{2\pi }}\int\limits_{B}{(\mathrm{div}\mathbf{A}_{a}-}\frac{{%
i\omega }}{c}\varphi _{a})\exp (-i\omega t)d\omega =0$

$\Delta \mathbf{A}-\frac{1}{{c^{2}}}\frac{{\partial ^{2}\mathbf{A}}}{{%
\partial t^{2}}}=\frac{1}{\sqrt{2\pi }}\int\limits_{D}{(\Delta \mathbf{A}%
_{a}+\frac{{\omega ^{2}}}{{c^{2}}}}\mathbf{A}_{a})\exp (-i\omega t)d\omega =-%
\frac{{4\pi }}{{\sqrt{2\pi }c}}\int\limits_{D}{\mathbf{j}_{a}}\exp (-i\omega
t)d\omega =-\frac{{4\pi }}{c}\mathbf{j}$

$\Delta \varphi -\frac{1}{{c^{2}}}\frac{{\partial ^{2}\varphi }}{{\partial
t^{2}}}=\frac{1}{\sqrt{2\pi }}\int\limits_{D}{(\Delta \varphi _{a}+\frac{{%
\omega ^{2}}}{{c^{2}}}}\varphi _{a})\exp (-i\omega t)d\omega =-\frac{{4\pi }%
}{\sqrt{2\pi }}\int\limits_{D}{\varphi _{a}}\exp (-i\omega t)d\omega =-4\pi
\rho $

\bigskip

i.e. (8), (9) and (10) have been proved.

The conditions of finitely extended Fourier transform on $\rho $ and $%
\mathbf{j}$ have been imposed since they are the simplest way to justify the
last passages of interchanging integrations and performing derivatives under
integration. However, the same justifications are ensured more generally by
conditions of existence and uniform absolute convergence of multiple
integrals of the form

$\bigskip $

$\int\limits_{B}{\int\limits_{D}{\frac{{\rho _{a}(r,\omega )}}{d}\omega ^{2}}%
}d\omega dV$

$\int\limits_{B}{\int\limits_{D}{\frac{{\mathbf{j}_{a}(r,\omega )}}{d}\omega
^{2}}}d\omega dV$

\bigskip

having considered that the factor $\exp (-i\omega t)$, of modulus 1, does
not disturb such absolute convergence.

\section{Conclusions}

A proof has been given that the fields calculated applying (6) and (7) to
the delayed potentials (14) and (15) provide effectively a solution to the
Maxwell equations (1), (2), (3), (4). The complete details of the proof have
been reported using the assumption that the sources $\rho $ and $\mathbf{j}$
are regular functions up to the first order derivatives, that they are
vanishing outside some finite spatial domain, and that their temporal
Fourier transform is vanishing outside some finite frequency interval. It
has also been shown that the last two conditions can be lowered, assuming
only the multiple summability of certain space/frequency functions
constructed from the sources.

\begin{flushleft}
\textbf{AMS Subject Classification: 78A25.}\\[2ex]
\end{flushleft}


\begin{thebibliography}{9}
\bibitem{}
 
 \textsc{J. D. Jackson}, \textit{Classical Electrodynamics}, J. Wiley and
Sons, N. Y. 1962.
\bibitem{}
\textsc{E. Durand}, \textit{Electrostatique et Magnetostatique}, Masson,
Paris 1953.
\bibitem{}
\textsc{J. H. Jeans}, \textit{Mathematical Theory of Electricity and
Magnetism}, Cambridgr University Press 1948.
\bibitem{}
\textsc{L. D. Landau}, \textit{Classical Theory of Fields},
Addison-Wesley 1951.
\bibitem{}
\textsc{E. Goursat}, \textit{Course d'Analyse Mathematique},
Gauthier-Villars, Paris 1949.
\end{thebibliography}
\end{document}